%% file: Dir_final_nod_arxiv.tex
\documentclass{elsart}

\usepackage{amssymb}

\newcommand\cal {}

\input macroarx

\begin{document}

\begin{frontmatter}

\title{Uniform asymptotic formulae for
Green's functions \\ in singularly perturbed domains} ~

\author{V. Maz'ya$^1$ and A. Movchan$^2$}

\address{$^1$ Department of Mathematical Sciences, University of Liverpool, \\ Liverpool L69 3BX,
U.K., ~ and
 Department of Mathematics, \\ Ohio State University, 231 W 18th Avenue, Columbus, \\ OH 43210,
 USA, ~ and
   Department of Mathematics, \\ Link\"oping University,  SE-581 83 Link\"oping, Sweden \\
~\\
  $^2$ Department of Mathematical Sciences, University of Liverpool, \\ Liverpool L69 3BX, U.K. \\
  ~ \\
  ~ \\
  {\rm Dedicated to Professor W.D. Evans } \\
  {\rm on the occasion of his sixty fifth birthday}}

\begin{abstract}
Asymptotic formulae for Green's functions for the operator $-\GD$
in
domains with small holes are obtained.
A new feature of these formulae is their uniformity with respect
to the independent variables.  The cases of multi-dimensional and
planar domains are considered.
\end{abstract}

\end{frontmatter}

%
%
%
%
%

\section{Introduction}

Hadamard's paper \cite{Had} contains, among much else, asymptotic
formulae for Green's kernels of classical boundary value problems
under  small variations of a domain. In \cite{Had}, the perturbed
domain $\GO_\Gve$, depending on a small parameter $\Gve > 0$,
approximates the limit domain $\GO$ in such a way that the angle
between the two outward normals at nearby
points of $\prt \GO$ and $\prt \GO_\Gve$ is small. In short,
Hadamard's formulae are related to the case of {\it a regularly
perturbed} domain. A drawback of these formulae is their
non-uniformity with respect to the independent variables. A
uniform version of one of Hadamard's formulae containing a
boundary layer was formulated in \cite{CRAS}. Besides, uniform
asymptotic representations of Green's functions for several types
of {\it singularly perturbed} domains were given in \cite{CRAS}
without proofs.

The objective of the present article is to prove two theorems
announced in \cite{CRAS}. Namely, we derive uniform asymptotic
formulae for Green's functions of the Dirichlet problem for the
operator $-\GD$ in $n$-dimensional domains with small holes, first
for $n>2$ in Section 2 and then for $n=2$ in Section 3.
Corollaries, presented in Section 4, show that these
formulae can be simplified under certain constraints on the
independent variables.


We make use of the version of the method of compound asymptotic
expansions of solutions to boundary value problems in singularly
perturbed domains
developed in
\cite{MNP}.
%
%
%

Now, we list several notations adopted in the text of the paper.
Let $\GO$ be a domain in ${\Bbb R}^n, ~ n \geq 2,$ with compact
closure $\ov{\GO}$ and boundary $\prt \GO$.  By $F$ we denote a
compact set of positive harmonic capacity in ${\Bbb R}^n$; its
complement is $F^c= {\Bbb R}^n \setminus F$. We suppose that both
$\GO$ and $F$ contain the origin $\BO$ as an interior point.
Without loss of generality, it is assumed that the minimum
distance between $\BO$ and  the points of $\prt \GO$ is equal to
$1.$ Also, the maximum distance between $\BO$ and the points of
$\prt F^c$ will be taken as $1$.
We introduce the set $F_\Gve = \{\Bx: \Gve^{-1} \Bx \in F \},$
where $\Gve$ is a small positive parameter, and
the open set $\GO_\Gve = \GO \setminus F_\Gve$. The notation
$B_\rho$ stands for  the open ball centered at $\BO$ with radius
$\rho$.
%
%
%
%

The main object of our study,  Green's function
for the operator $-\GD$ in $\GO_\Gve$, will be denoted by
$G_\Gve$.
%
%
In the sequel, along with $\Bx$ and $\By$, we use the  scaled
variables $\BGx = \Gve^{-1} \Bx~~ \mbox{and} ~~ \BGn = \Gve^{-1}
\By.$

By $\mbox{Const}$ we always mean different positive constants
depending only on $n$. Finally, the notation $f=O(g)$ is
equivalent to the inequality $|f| \leq \mbox{Const} ~ g$.

\section{Green's function for a multi-dimensional domain with a small hole}
\label{s1}

%
%

We assume here that $n >2.$ Let $G$ and $g$ denote Green's
functions of the Dirichlet problem for the operator $-\GD$ in the
sets $\GO$ and $F^c={\Bbb R}^n \setminus F$.
We make use of the regular parts of $G$ and $g$, respectively:
\beq {\cal H}(\Bx, \By)  = {(n-2)^{-1} |S^{n-1}|^{-1} |\Bx -
\By|^{2-n}} - G(\Bx, \By), \eequ{eq3} and \beq h(\BGx, \BGn) =
{(n-2)^{-1} |S^{n-1}|^{-1} |\BGx - \BGn|^{2-n}} - g(\BGx, \BGn),
\eequ{eq3b} where $|S^{n-1}|$ denotes the $(n-1)$-dimensional
measure of the unit sphere $S^{n-1}$.


By $P(\BGx)$ we mean the equilibrium potential  of $F$ defined as
a unique solution of the following Dirichlet problem in
$F^c$ \beq \GD_\Gx P(\BGx) = 0 ~~\mbox{in} ~  F^c,
\eequ{eq7} \beq P(\BGx) = 1 ~~ \mbox{on} ~ \prt F^c,
\eequ{eq8} \beq P(\BGx) \to 0 ~~\mbox{as} ~~ |\BGx|
\to \infty, \eequ{eq9} where the boundary condition \eq{eq8} is
interpreted in the sense of the Sobolev space $H^1$.

%

%

%
%
%


The following auxiliary assertion is classical.

{\bf Lemma 1.}

{\em {\it (i)}  The potential $P$ satisfies the estimate \beq 0 <
P(\BGx) \leq \min \Big\{ 1,
{|\BGx|}
^{2-n} \Big\}. \eequ{P_eq1}

{\it (ii)} If $|\BGx| \geq 2$,
then \beq \Big|P(\BGx) -
\fr{\mbox{\rm cap}(F)}{(n-2) |S^{n-1}|} |\BGx|^{2-n} \Big| \leq
\mbox{\rm Const} ~
{|\BGx|}^{1-n}
\eequ{P_eq2}

}

{\bf Proof.} {\it (i)} Inequalities \eq{P_eq1} follow from the
maximum principle for variational solutions of Laplace's equation.

{\it (ii)} Inequality \eq{P_eq2} results from the expansion of $P$
in spherical harmonics. $\Box$

{\bf Lemma 2.} {\em For all $\BGn \in F^c$ and for $\BGx$ with
$|\BGx| > 2$
the estimate holds:}
\beq |h(\BGx, \BGn) - P(\BGn){(n-2)^{-1} |S^{n-1}|^{-1}
|\BGx|^{2-n}}| \leq
\mbox{\rm Const}~
|\BGx|^{1-n} P(\BGn).
\eequ{eq11}

{\bf Proof.} By {\rm \eq{eq3b}}, $h$ satisfies the Dirichlet
problem \beqa \GD_\Gx h(\BGx, \BGn) &=& 0, ~~ \BGx, \BGn \in
F^c,
\label{e1} \\
 h(\BGx, \BGn) &=& (n-2)^{-1} |S^{n-1}|^{-1} |\BGx - \BGn|^{2-n}, \nonumber \\
&& ~~ \BGx \in \prt F^c
~\mbox{and} ~ \BGn \in
F^c,
\label{e2} \\
h(\BGx, \BGn) &\to& 0 ~~\mbox{as} ~ |\BGx| \to \infty ~\mbox{and}~
\BGn \in
F^c. \label{e3}
\end{eqnarray}


We fix $\BGn \in F^c.$ By the series expansion of $g$ in spherical
harmonics,
\beq |\BGx|^{n-2} \Big(g(\BGx, \BGn) - \fr{C(\BGn)}{(n-2)
|S^{n-1}| |\BGx|^{n-2}}\Big) \to 0 ~~
\mbox{as} ~ |\BGx| \to \infty. \eequ{g}



We apply Green's formula to the functions $g(\BGx, \BGn)$ and $1 -
P(\BGx)$ restricted to the domain $B_R \setminus F$, where $B_R =
\{ \BGx: |\BGx| < R \}$ is the ball of a sufficiently large radius
$R$. Taking into account that $P(\BGx) =1$ and $g(\BGx, \BGn)=0$
when  $\BGx \in\prt (F^c)$ we deduce \beq \hspace{-.2in} \int_{B_R
\setminus F} \Grad_\Gx g(\BGx, \BGn) \cdot \Grad_\Gx P(\BGx)  d
\BGx = P(\BGn)-1 - \int_{\prt B_R} (1-P(\BGx)) \fr{\prt}{\prt
|\BGx|}g(\BGx, \BGn) d s_\Gx, \eequ{h1} and \beq \int_{B_R
\setminus F} \Grad_\Gx g(\BGx, \BGn) \cdot \Grad_\Gx P(\BGx)  d
\BGx = \int_{\prt B_R} g(\BGx, \BGn) \fr{\prt}{\prt |\BGx|}
P(\BGx) d s_\Gx. \eequ{h2} Hence, \beq 1-P(\BGn) = -\int_{\prt
B_R} \bigg( g(\BGx, \BGn) \fr{\prt}{\prt |\BGx|} P(\BGx)+
(1-P(\BGx)) \fr{\prt}{\prt |\BGx|}g(\BGx, \BGn) \bigg) d s_\Gx.
\eequ{h3} It follows from
\eq{g}
that
$$
1-P(\BGn) = - \lim_{R \to \infty} \int_{\prt B_R} \fr{\prt}{\prt
|\BGx|} \fr{C(\BGn)}{(n-2) |S^{n-1}| |\BGx|^{n-2}} d s_\Gx =
C(\BGn).
$$

%
%
%

Let
$|\BGx| > 2$.
Then
for $\BGn \in \prt F^c$
$$
|h(\BGx, \BGn) - (n-2)^{-1} |S^{n-1}|^{-1} |\BGx|^{2-n} P(\BGn) |
=(n-2)^{-1} |S^{n-2}|^{-1} \Big||\BGx-\BGn|^{2-n} -
|\BGx|^{2-n}\Big|
$$
\beq \leq \mbox{Const} ~ |\BGn| |\BGx|^{1-n} \leq \mbox{Const} ~
|\BGx|^{1-n}. \eequ{h_est_inf} In the above estimate, we used the
assumption (see Introduction) of the maximum distance between the
origin and the points of $\prt F^c$ being equal to $1$. From
\eq{h_est_inf} and the maximum principle for functions harmonic in
$\BGn$, we deduce
$$
|h(\BGx, \BGn) - \Big((n-2)|S^{n-1}|\Big)^{-1}|\BGx|^{2-n}
P(\BGn)| \leq \mbox{Const}~  |\BGx|^{1-n} P(\BGn),
$$
for all $\BGn \in F^c$ and $|\BGx| > 2 $. $\Box$

Our main result concerning the uniform approximation of Green's
function $G_\Gve$ in the multi-dimensional case is given by

{\bf Theorem 1.} {\em Green's function $G_\Gve(\Bx, \By)$
admits the
representation
\begin{eqnarray}
G_\Gve(\Bx, \By) &&= G(\Bx, \By) + \Gve^{2-n}g(\Gve^{-1}\Bx,
\Gve^{-1}\By) - ((n-2) |S^{n-1}| |\Bx - \By|^{n-2})^{-1}
\nonumber  \\
  +{\cal H}(0, &&\By) P(\Gve^{-1} \Bx) + {\cal H}(\Bx, 0)
P(\Gve^{-1} \By)
- {\cal H}(0,0) P(\Gve^{-1} \Bx) P(\Gve^{-1} \By)
\nonumber \\
 -\Gve^{n-2} &&\mbox{\rm cap}(F) ~{\cal H}(\Bx, 0) {\cal H}(0, \By)
+ O\Big(
{\Gve^{n-1}}{(\min \{|\Bx|, |\By| \}+ \Gve)^{2-n}}\Big),
\label{eq16}
\end{eqnarray}
uniformly with respect to
$\Bx, \By \in \GO_\Gve.$
Here, ${\cal H}$ and $h$ are regular parts of Green's functions
$G$ and $g$, respectively (see {\rm \eq{eq3}, \eq{eq3b}}), and $P$
is the capacitary potential of $F$.}


Before presenting a proof of this theorem, we give a {\it
plausible formal argument} leading to  \eq{eq16}.

Let $G_\Gve$
be represented in the form \beq G_\Gve(\Bx, \By) =
\Big((n-2) |S^{n-1}|\Big)^{-1} |\Bx - \By|^{2-n}
- {\cal H}_\Gve(\Bx,
\By) -
h_\Gve(\Bx, \By), \eequ{eq6} where ${\cal H}_\Gve$ and $h_\Gve$
are solutions of the Dirichlet
problems
\beqa
\GD_x {\cal H}_\Gve(\Bx, \By) = 0, ~~ \Bx, \By \in \GO_\Gve, \nonumber \\
{\cal H}_\Gve(\Bx, \By) = \Big((n-2) |S^{n-1}|\Big)^{-1} |\Bx -
\By|^{2-n}, ~~ \Bx \in \prt \GO, ~\By \in \GO_\Gve,
\nonumber \\
{\cal H}_\Gve(\Bx, \By) = 0,  ~~ \Bx \in \prt F^c_\Gve,
~\By \in \GO_\Gve. \nonumber
\end{eqnarray}
and
\beqa
\GD_x h_\Gve(\Bx, \By) = 0,  ~~ \Bx, \By \in \GO_\Gve, \nonumber \\
h_\Gve(\Bx, \By) =
\Big((n-2) |S^{n-1}| \Big)^{-1} |\Bx - \By|^{2-n}, ~~ \Bx \in \prt
F^c_\Gve,
~\By \in \GO_\Gve,
\label{eqh_ep_d} \\
h_\Gve(\Bx, \By) = 0,  ~~ \Bx \in \prt \GO, ~\By \in \GO_\Gve.
\nonumber
\end{eqnarray}




By \eq{eq6}, it suffices to find asymptotic formulae for ${\cal
H}_\Gve$ and $h_\Gve$.

{\em Function ${\cal H}_\Gve$.}
Obviously, ${\cal H}_\Gve(\Bx, \By)-{\cal H}(\Bx, \By)$ is
harmonic in $\GO_\Gve$, and ${\cal H}_\Gve(\Bx, \By)-{\cal H}(\Bx,
\By) =0$ for $\Bx \in \prt \GO$.
On the other hand, for $\Bx \in \prt F^c_\Gve$ the leading part of
${\cal H}_\Gve(\Bx, \By)-{\cal H}(\Bx, \By)$
is equal to the function $-{\cal H}(0, \By)$.
This function can be  extended onto $F^c_\Gve$, harmonically in
$\Bx$, as $-{\cal H}(0, \By) P(\Gve^{-1} \Bx)$, whose
leading-order part is equal to $-\Gve^{n-2} \mbox{cap}(F) ~ {\cal
H}(\Bx, 0) {\cal H}(0, \By)$ for $\Bx \in \prt \GO$.
 Hence,
$$
{\cal H}_\Gve(\Bx, \By) - {\cal H}(\Bx, \By) \sim - {\cal H}(0,
\By) P(\Gve^{-1} \Bx)
$$
\beq ~~~~~~~~~+ \Gve^{n-2} \mbox{cap}(F) ~ {\cal H}(\Bx, 0) {\cal
H}(0, \By) ~~ \mbox{for all}~ \Bx, \By \in \GO_\Gve. \eequ{eq12na}

{\em Function $h_\Gve$.} By definitions \eq{eq3b} and
\eq{eqh_ep_d} of $h$ and $h_\Gve$,
$$
h_\Gve(\Bx, \By)-\Gve^{2-n} h(\Gve^{-1}\Bx, \Gve^{-1}\By)=0
~\mbox{for} ~ \Bx \in \prt F^c_\Gve.
$$
Furthermore, by Lemma 2
$$
h_\Gve(\Bx, \By) -
\Gve^{2-n} h(\Gve^{-1} \Bx, \Gve^{-1} \By)
$$
$$
\sim - \Big((n-2) |S^{n-1}| \Big)^{-1} |\Bx|^{2-n} P(\Gve^{-1}
\By) ~~ \mbox{for}~ \Bx \in \prt \GO.
$$
The harmonic function in $\Bx \in \GO$, with the Dirichlet data
%
%
$$- \Big((n-2) |S^{n-1}| \Big)^{-1} |\Bx|^{2-n} P(\Gve^{-1} \By) $$ on $\prt
\GO$, is $-{\cal H}(\Bx, 0) P(\Gve^{-1} \By)$, and it is
asymptotically equal to $-{\cal H}(0,0) P(\Gve^{-1} \By)$ on $\prt
F^c_\Gve$,
which is not necessarily small.
%
The harmonic in $\Bx$ extension  of  \\ ${\cal H}(0,0) P(\Gve^{-1}
\By)$ onto $F^c_\Gve$ is given by ${\cal H}(0,0) P(\Gve^{-1} \By)
P(\Gve^{-1} \Bx)$. Since this function is small for $\Bx \in \prt
\GO$, one may assume the asymptotic representation
%
%
%
%
\begin{eqnarray}
h_\Gve(\Bx, \By) &&- \Gve^{2-n} h(\Gve^{-1} \Bx, \Gve^{-1} \By) +
{\cal H}(\Bx, 0)
P(\Gve^{-1} \By) \nonumber \\
&& \sim
{\cal H}(0,0) P(\Gve^{-1} \Bx) P(\Gve^{-1} \By)~~~~~~~ \mbox{for
all}~ \Bx, \By \in \GO_\Gve.
\label{eq12m}
\end{eqnarray}


Substituting \eq{eq12na} and \eq{eq12m} into
\eq{eq6}, we deduce
\begin{eqnarray}
G_\Gve(\Bx, \By) && \sim \Big((n-2) |S^{n-1}| \Big)^{-1} |\Bx -
\By|^{2-n} -{\cal H}(\Bx, \By) - \Gve^{2-n} h(\Gve^{-1}\Bx,
\Gve^{-1}\By)
\nonumber \\
&& +{\cal H}(0, \By) P(\Gve^{-1} \Bx) + {\cal H}(\Bx, 0)
P(\Gve^{-1} \By) - {\cal H}(0,0) P(\Gve^{-1} \Bx) P(\Gve^{-1} \By)
\nonumber \\
&& -\Gve^{n-2} \mbox{cap}(F) ~ {\cal H}(\Bx, 0) {\cal H}(0, \By),
\nonumber
\end{eqnarray}
which is equivalent to
\begin{eqnarray}
G_\Gve(\Bx, \By) &&\sim G(\Bx, \By) + \Gve^{2-n}g(\Gve^{-1}\Bx,
\Gve^{-1}\By) - ((n-2) |S^{n-1}|)^{-1} |\Bx - \By|^{2-n}
\nonumber \\
&& +{\cal H}(0, \By) P(\Gve^{-1} \Bx) + {\cal H}(\Bx, 0)
P(\Gve^{-1} \By)
- {\cal H}(0,0) P(\Gve^{-1} \Bx) P(\Gve^{-1} \By)
\nonumber \\
 && -\Gve^{n-2} \mbox{cap}(F) ~ {\cal H}(\Bx, 0) {\cal
H}(0, \By). \nonumber
\end{eqnarray}
Now, we give a rigorous proof of \eq{eq16}.



{\bf Proof of Theorem 1.}

The remainder $r_\Gve(\Bx, \By)$ in \eq{eq16}
is a solution of the boundary value problem
\beqa
\GD_x r_\Gve(\Bx, \By) &=& 0, ~~ \Bx, \By \in \GO_\Gve,
\label{eq17} \\
\nonumber \\
r_\Gve(\Bx, \By) &=& {\cal H}(\Bx, \By) - {\cal H}(0, \By) \nonumber \\
&& -({\cal H}(\Bx,0)-{\cal H}(0,0)) P(\Gve^{-1} \By) \nonumber \\
&& + \Gve^{n-2} \mbox{\rm cap}(F) ~ {\cal H}(\Bx, 0) {\cal H}(0,
\By), ~~ \Bx \in \prt
F_\Gve^c, ~ \By \in \GO_\Gve,
\label{eq19} \\
\nonumber \\
 r_\Gve(\Bx, \By) &=& \Gve^{2-n} h(\Gve^{-1} \Bx,
\Gve^{-1} \By) - {\cal H}(0, \By) P(\Gve^{-1}\Bx) \nonumber
\\
&-&{\cal H}(\Bx, 0) P(\Gve^{-1} \By) + {\cal H}(0,0) P(\Gve^{-1}
\Bx) P(\Gve^{-1} \By)
\nonumber \\
&& + \Gve^{n-2} \mbox{\rm cap}(F) ~ {\cal H}(\Bx, 0) {\cal H}(0,
\By), ~ \Bx \in \prt \GO, ~ \By \in \GO_\Gve. \label{eq18}
\end{eqnarray}

The functions ${\cal H}(\Bx,0)$ and ${\cal H}(0, \By)$ are
harmonic in $\GO$ and are
bounded by $\mbox{Const}$
on $\prt \GO$. Hence, they are bounded by $\mbox{Const}$ for $\Bx
\in \prt F_\Gve^c, ~ \By \in \GO_\Gve$ and for $\Bx \in \prt \GO,
~ \By \in \GO_\Gve,$ respectively. The terms $\Gve^{n-2} \mbox{\rm
cap}(F){\cal H}(\Bx, 0) {\cal H}(0, \By)$ in the right-hand sides
of \eq{eq19} and \eq{eq18} are
bounded by $\mbox{Const}~ \Gve^{n-2}$.

By definition \eq{eq3}, $\Grad_x {\cal H}(\Bx, \By)$ is bounded by
$\mbox{Const}$ uniformly with respect to $\By \in \GO$ for every
$\Bx \in B_{1/2}$.
Hence, by \eq{eq19} and the inequalities $0 < P(\Bx) \leq 1$,
$$
|{\cal H}(\Bx, \By) - {\cal H}(0, \By) -({\cal H}(\Bx,0)-{\cal
H}(0,0)) P(\Gve^{-1} \By)|
$$
$$
\leq  \mbox{Const}~\Gve  ~\sup_{\Bz \in B_{\Gve}} |\Grad_z {\cal
H}(\Bz, \By)| \leq \mbox{Const} ~ \Gve ,
$$
for $\Bx \in \prt F^c_\Gve, ~ \By \in \GO_\Gve$.
Thus, the following estimate holds
when $\Bx \in \prt
F^c_\Gve$ and $\By \in \GO_\Gve$ \beq |r_\Gve (\Bx, \By)| \leq
\mbox{Const} ~ \Gve  \sup_{\Bz \in B_{\Gve }} |\Grad_z {\cal
H}(\Bz, \By)| \leq \mbox{Const} ~ \Gve . \eequ{F_est}
Next, we estimate
 $|r_\Gve(\Bx, \By)|$ for $\Bx \in \prt \GO$ and $\By
\in \GO_\Gve$.
By Lemma 1, the capacitary potential $P(\Gve^{-1} \Bx)$ satisfies
the inequalities \beq 0 \leq P(\Gve^{-1} \Bx) \leq \mbox{\rm
Const} ~\fr{\Gve^{n-2}}{(|\Bx| + \Gve )^{n-2}}, \eequ{P_est}
for $\Bx \in \GO_\Gve$, and \beqa &&\Big|P(\Gve^{-1} \Bx) -
\fr{\Gve^{n-2} \mbox{cap}( F)}{(n-2)
|S^{n-1}| |\Bx|^{n-2}} \Big| \nonumber \\
&& \leq \mbox{\rm Const} ~ \Big( \Gve /{|\Bx|} \Big)^{n-1} \leq
\mbox{\rm Const} ~ \Gve^{n-1}, \label{Cap_pot}
\end{eqnarray}
for $\Bx \in \prt \GO$. Now, \eq{Cap_pot} and  the definition of
$H(\Bx, \By)$ imply \beq |\Gve^{n-2} \mbox{\rm cap}(F) H(\Bx,0)
H(0, \By) - H(0,\By) P(\Gve^{-1} \Bx)| \leq \mbox{\rm Const} ~
\Gve^{n-1}. \eequ{Cap_pot1}

Also, we have the estimate \beqa && |\Gve^{2-n} h(\Gve^{-1} \Bx,
\Gve^{-1} \By) - {\cal H}(\Bx, 0) P(\Gve^{-1} \By) | \nonumber \\
&&=\Gve^{2-n} \Big| h(\Gve^{-1} \Bx, \Gve^{-1} \By) -
\fr{P(\Gve^{-1} \By)}{(n-2) |S^{n-1}| |\Bx / \Gve|^{n-2}} \Big| \nonumber \\
&& \leq \mbox{\rm Const} ~\Gve  |\Bx|^{1-n} P(\Gve^{-1} \By) \nonumber \\
&& \leq \mbox{Const} ~
\fr{\Gve^{n-1}}{(|\By| + \Gve)^{n-2}},
~ \Bx \in \prt \GO, ~ \By \in \GO_\Gve, \label{hh_est}
\end{eqnarray}
which follows from the definition \eq{eq3} of $H(\Bx, \By)$  and
the estimates \eq{eq11} and \eq{P_est}.
Combining \eq{P_est}, \eq{Cap_pot1} and \eq{hh_est}  we obtain
from \eq{eq18} that the trace of the function $\Bx \to
|r_\Gve(\Bx, \By)|$ on $\prt \GO$ does not exceed $$
\mbox{\rm Const} ~ \fr{\Gve^{n-1}}{(|\By| + \Gve )^{n-2}}.
$$
for $
\By \in \GO_\Gve$.
Using this and \eq{F_est}, we deduce by the maximum principle
that
$$
|r_\Gve(\Bx, \By)| \leq \mbox{\rm Const} \Big\{ \Gve
P\Big(\fr{\Bx}{\Gve}\Big) + \fr{\Gve^{n-1}}{(|\By| + \Gve)^{n-2}}
\Big\},
$$
for all $\Bx, \By \in \GO_\Gve$. Taking into account \eq{P_est},
we arrive at \beq |r_\Gve(\Bx, \By)| \leq \mbox{\rm Const}~
\fr{\Gve^{n-1}}{(\min \{|\Bx|, |\By| \}+ \Gve)^{n-2}}
\eequ{est_est}


The proof is complete. $\Box$

\section{Green's function for the Dirichlet problem in a planar domain
with a small hole}
\label{s3}

In this section, we find  an asymptotic approximation of $G_\Gve$
in the two-dimensional case. We shall see that this approximation
has new features in comparison with that in Theorem 1.

The notations $\GO_\Gve, \GO, F_\Gve, F,$ introduced in
Introduction, will be used here. As before, we assume that the
minimum distance from the origin to $\prt \GO$
and the maximum distance between the origin and the points of
$\prt F^c$ are equal to $1$.

%

Green's function $G(\Bx,
\By)$ for the unperturbed domain $\GO$  has the form \beq G(\Bx,
\By) =
{(2 \pi)^{-1}} \log |\Bx - \By|^{-1} - H(\Bx, \By), \eequ{eq32}
where $H$ is its regular part satisfying
\beq
\GD_x H (\Bx, \By) = 0,  ~~ \Bx, \By \in \GO,
\eequ{eq33}
\beq
H(\Bx, \By) =
{(2 \pi)^{-1}} \log |\Bx - \By|^{-1}, ~~ \Bx \in \prt \GO, ~\By
\in \GO. \eequ{eq34}

The scaled coordinates $\BGx=\Gve^{-1} \Bx$ and $\BGn= \Gve^{-1}
\By$ will be used
as in the multi-dimensional case. Similar to Section \ref{s1},
$g(\BGx, \BGn)$ and
$h(\BGx, \BGn)$
are Green's function and its regular part in $F^c $:
\beqa \GD_\Gx g (\BGx, \BGn) + \Gd(\BGx - \BGn) = 0, ~~ \BGx, \BGn
\in F^c,
\label{eq38a} \\
\nonumber \\
 g(\BGx, \BGn) = 0, ~~ \BGx \in \prt F,~
\BGn \in F^c,
\label{eq39a} \\
\nonumber \\
g(\BGx, \BGn) ~~\mbox{is bounded} ~~\mbox{as} ~ |\BGx| \to \infty
~\mbox{and}~ \BGn \in F^c,
\label{eq40a}
\end{eqnarray}
and
\beq
h(\BGx, \BGn) =
{(2 \pi)^{-1}} \log |\BGx - \BGn|^{-1} - g(\BGx, \BGn).
\eequ{eq38}

We introduce a function $\Gz$
by
\beq
\Gz(\BGn) = \lim_{|\BGx|\to \infty} g(\BGx, \BGn),
\eequ{eq40aa}
and the constant
\beq
\Gz_\infty = \lim_{|\BGn| \to \infty} \{ \Gz(\BGn) -
{(2 \pi)^{-1}} \log |\BGn| \}. \eequ{eq40b}


{\bf Lemma 3.} {\em Let $|\BGx| > 2$. Then the  regular part
$h(\BGx, \BGn)$ of Green's function $g$ in $F^c$
admits the asymptotic representation \beq h(\BGx, \BGn) = -
{(2 \pi)^{-1}} \log |\BGx| - \Gz(\BGn) + O(|\BGx|^{-1}),
\eequ{h_eq} which is uniform with respect to $\BGn \in F^c$.}

{\bf Proof:} Following the inversion transformation, we  use the
variables:
$$
\BGx' = |\BGx|^{-2} \BGx, ~~ \BGn' = |\BGn|^{-2} \BGn,
$$
and the identity
$$
|\BGx - \BGn|^{-1}|\BGx||\BGn| = |\BGx'-\BGn'|^{-1}.
$$
Then, the boundary values of $h(\BGx, \BGn)$, as $\BGx \in \prt
F^c, \BGn \in F^c$,
can be expressed in the form \beq ~~~~~~~~~~~~~~h(\BGx, \BGn) =
{\frak H}(\BGx', \BGn') -
{(2 \pi)^{-1}} \log |\BGx||\BGn|, \eequ{bv_h} where ${\frak
H}(\BGx', \BGn'), ~\BGx' \in \prt (F^c)',$
is the boundary value of the regular part of Green's function in
the bounded transformed
set $(F^c)'$.
Namely, the function ${\frak H}(\BGx', \BGn')$ is defined as a
solution of the Dirichlet problem \beq ~~~~~~~~~~~~ \GD_{\BGx'}
{\frak H}(\BGx', \BGn') = 0, ~~ \BGx', \BGn' \in (F^c)',
\eequ{hi1} \beq ~~~~~~~~~~~~ {\frak H}(\BGx', \BGn') =
{(2 \pi)^{-1}} \log |\BGx' - \BGn'|^{-1}, ~~ \BGx' \in \prt
(F^c)'.
\eequ{hi2}

It follows from \eq{bv_h} that the harmonic extension of $h(\BGx,
\BGn)$ is \beq h(\BGx, \BGn) = {\frak H}(\BGx', \BGn') -
{(2 \pi)^{-1}} \log |\BGx||\BGn|, ~~ \BGx, \BGn \in F^c.
\eequ{h_FF} Since ${\frak H}(\BGx', \BGn')$ is smooth in $(F^c)'
\times (F^c)'$, we deduce \beq h(\BGx, \BGn) = {\frak H}(0, \BGn')
-
{(2 \pi)^{-1}} \log|\BGx||\BGn| + O(|\BGx'|),
\eequ{h_asym}
for $|\BGx'| < 1/2$ and for all $\BGn' \in (F^c)'.$ Also, by
\eq{h_FF} and the definition of $h(\BGx, \BGn)$, \beq {\frak
H}(\BGx', \BGn')= -g(\BGx, \BGn) + (2 \pi)^{-1} \log |\BGx| |\BGn|
- (2 \pi)^{-1} \log |\BGx - \BGn|. \eequ{g_h} Then, applying
\eq{eq40aa} and taking the limit in \eq{g_h}, as $|\BGx'| \to 0$,
we arrive at
$${\frak H}(0, \BGn')= - \Gz(\BGn) + (2 \pi)^{-1}\lim_{|\BGx| \to
\infty} \log (|\BGx - \BGn|^{-1} |\BGx|) + (2 \pi)^{-1} \log
|\BGn|
$$
$$
={(2 \pi)^{-1}} \log |\BGn| - \Gz(\BGn).$$ Further substitution of
${\frak H}(0, \BGn')$ into \eq{h_asym} leads to
$$
h(\BGx, \BGn) = -
{(2 \pi)^{-1}} \log |\BGx| - \Gz(\BGn) + O(|\BGx|^{-1}),
$$
for $|\BGx| > 2$ and for all $\BGn \in F^c.$ The proof is complete
$\Box$.

\subsection{Asymptotic approximation of the equilibrium potential}
\label{s3a}

The {\em equilibrium potential} $\CP_\Gve(\Bx)$ is introduced as a
solution of
the following Dirichlet problem in $\GO_\Gve$
\beqa
\GD \CP_\Gve (\Bx) = 0, ~~ \Bx \in \GO_\Gve, \label{eq41} \\
\CP_\Gve(\Bx) = 0, ~~\Bx \in \prt \GO, \label{eq42} \\
\CP_\Gve(\Bx) = 1, ~~\Bx \in \prt F_\Gve^c. \label{eq43}
\end{eqnarray}

{\bf Lemma 4.} {\em The asymptotic approximation of
$\CP_\Gve(\Bx)$ is given by the formula \beq \CP_\Gve (\Bx) =
\fr{-G(\Bx, 0) + \Gz(\fr{\Bx}{\Gve}) - \fr{1}{2 \pi} \log
\fr{|\Bx|}{\Gve} - \Gz_\infty}{\fr{1}{2 \pi} \log \Gve + H(0, 0) -
\Gz_\infty} +p_\Gve(\Bx), \eequ{eq44} where
$\Gz_\infty$ is
defined by {\rm \eq{eq40b}}, and $p_\Gve$ is the remainder term
such that
$$
|p_\Gve(\Bx)| \leq \mbox{Const} ~ \Gve (\log \Gve)^{-1}
$$
uniformly with
respect to $\Bx \in \GO_\Gve$.}

{\bf Proof.} Direct substitution of \eq{eq44} into
\eq{eq41}--\eq{eq43} yields the Dirichlet problem for the
remainder term $p_\Gve$ \beqa
\GD p_\Gve(\Bx) &=& 0, ~~\Bx \in \GO_\Gve, \label{eq45} \\
p_\Gve(\Bx) &=& -\fr{\Gz(
{\Gve}^{-1} {\Bx})  - \fr{1}{2 \pi} \log
({\Gve}^{-1}{|\Bx|}) - \Gz_\infty}{\fr{1}{2 \pi} \log \Gve +
H(0,0) - \Gz_\infty},
%
~~\Bx \in \prt \GO, \label{eq46} \\
p_\Gve(\Bx) &=& 1-\fr{H(\Bx, 0) + \fr{1}{2 \pi} \log {\Gve} -
\Gz_\infty}{\fr{1}{2 \pi} \log \Gve + H(0,0) - \Gz_\infty},~~\Bx
\in \prt F_\Gve^c. \label{eq47}
\end{eqnarray}

Using \eq{eq40b} and the expansion of $\Gz(\BGx)$ in spherical
harmonics, we deduce
$$
\Gz(\Gve^{-1} \Bx) -
{(2 \pi)^{-1}} \log (\Gve^{-1} |\Bx|) - \Gz_\infty = O(\Gve ),
$$
as $|\Bx| \in \prt \GO,$ and hence the right-hand side in
\eq{eq46} is $O(\Gve(\log \Gve)^{-1})$. Since $H(\Bx, 0)$ is
smooth in $\GO$, we have
$$
H(\Bx,0) - H(0,0) = O(\Gve),
$$
as $\Bx \in \prt F^c_\Gve,$ and therefore
the right-hand side in \eq{eq47} is also $O(\Gve(\log
\Gve)^{-1})$. Applying the maximum principle, we arrive at the
result of Lemma. $\Box$

{\bf Remark.} For the case when
$\GO$ is a Jordan domain and $F$ is the closure of a Jordan
domain, we can adopt the notions of \cite{PS}: the inner conformal
radius $r_F$  of $F$, with respect to $\BO$, and the outer
conformal radius $R_\GO$ of $\GO$, with respect to $\BO$, are
defined as
$$
r_F = \exp (- 2 \pi \zeta_\infty), ~R_\GO = \exp (- 2 \pi H(0,0)),
$$
respectively.
%
In this case, the equilibrium potential $\CP_\Gve(\Bx)$ can be
represented in the form
$$
\CP_\Gve (\Bx) = \fr{-G(\Bx, 0) + \Gz(\fr{\Bx}{\Gve}) - \fr{1}{2 \pi} \log
\fr{|\Bx|}{\Gve r_F}}{\fr{1}{2 \pi} \log \fr{\Gve r_F}{R_\GO} }
+p_\Gve(\Bx).
$$

\subsection{Uniform asymptotic approximation}


{\bf Theorem 2.} {\em Green's function $G_\Gve$ for  the operator
$-\GD$
in $\GO_\Gve \subset {\Bbb R}^2$ admits the representation}
$$
G_\Gve(\Bx, \By) = G(\Bx, \By)  + {g}(\Gve^{-1} \Bx, \Gve^{-1}\By)
+ (2 \pi )^{-1} \log(\Gve^{-1} |\Bx - \By| )
$$
$$
+\fr{\Big( (2 \pi)^{-1}\log \Gve + \zeta(\fr{\Bx}{\Gve})
-\zeta_\infty+ H(\Bx, 0) \Big) \Big( (2 \pi)^{-1}\log \Gve +
\zeta(\fr{\By}{\Gve}) -\zeta_\infty + H(0, \By) \Big)}{(2
\pi)^{-1} \log \Gve + H(0,0) - \zeta_\infty}
$$
\beq - \zeta(\Gve^{-1} \Bx) - \zeta(\Gve^{-1} \By) + \zeta_\infty
+O(\Gve), \eequ{eq48}
%
{\em which is uniform with respect to $(\Bx, \By) \in \GO_\Gve
\times \GO_\Gve.$}

{\bf Proof.}
Let
\beq
G_\Gve(\Bx, \By) =
{(2 \pi)^{-1}} \log|\Bx-\By|^{-1} -H_\Gve(\Bx, \By) - h_\Gve(\Bx,
\By), \eequ{eq48a} where $H_\Gve$ and $h_\Gve$ are defined as
solutions of the Dirichlet problems \beqa \GD_x H_\Gve(\Bx, \By) =
0, ~~\Bx,~ \By \in \GO_\Gve,
\label{eq49} \\
H_\Gve(\Bx, \By) =
{(2 \pi)^{-1}} \log |\Bx-\By|^{-1}, ~~\Bx \in \prt \GO, ~ \By \in
\GO_\Gve,
\label{eq50} \\
H_\Gve(\Bx, \By) = 0, ~~ \Bx \in \prt F_\Gve, ~ \By \in \GO_\Gve,
\label{eq51}
\end{eqnarray}
and
 \beqa
\GD_x h_\Gve(\Bx, \By) = 0, ~~\Bx,~ \By \in \GO_\Gve,
\label{eq52} \\
h_\Gve(\Bx, \By) = 0,
~~\Bx \in \prt \GO, ~ \By \in \GO_\Gve,
\label{eq53} \\
h_\Gve(\Bx, \By) =
{(2 \pi)^{-1}} \log |\Bx-\By|^{-1}, ~~ \Bx \in \prt F_\Gve, ~ \By
\in \GO_\Gve. \label{eq54}
\end{eqnarray}

{\em The function $H_\Gve$} is represented in the form
\beq
H_\Gve(\Bx, \By) = C(\By, \log \Gve) G(\Bx, 0) + H(\Bx, \By) +
R_\Gve(\Bx, \By, \log \Gve),
\eequ{eq54a}
where
$C(\By, \log \Gve)$ is to be determined,
$G$ and $H$ are defined by \eq{eq32}--\eq{eq34}, and the third
term $R_\Gve$ satisfies the boundary value problem
\beqa
\GD_x R_\Gve(\Bx, \By, \log \Gve) = 0, ~~\Bx,~ \By \in \GO_\Gve, \label{eq55} \\
R_\Gve(\Bx, \By, \log \Gve) = 0, ~~ \Bx \in \prt \GO, ~\By \in \GO_\Gve, \label{eq56} \\
R_\Gve(\Bx, \By, \log \Gve) = - C G(\Bx, 0) - H(\Bx, \By), ~~ \Bx
\in \prt F_\Gve, ~\By \in \GO_\Gve, \label{eq57}
\end{eqnarray}
and it is approximated by a function $R(\Gve^{-1} \Bx, \By, \log
\Gve)$ defined in scaled coordinates in such a way that \beqa
\GD_\Gx R(\BGx,\By, \log \Gve) =&& 0, ~~ \BGx \in F^c,
\label{eq58} \\
R(\BGx, \By, \log \Gve) =&&
{C}{(2 \pi)^{-1}}( \log |\BGx| +
\log \Gve)
\nonumber \\
&& + C H(0,0) - H(0, \By),  ~~\BGx \in \prt F^c,
\label{eq59} \\
R(\BGx, \By, \log \Gve) \to && 0~~ \mbox{as} ~ |\BGx| \to \infty,
\label{eq60}
\end{eqnarray}
where $\By \in \GO_\Gve.$
The solution of the above problem has the form \beqa R(\BGx, \By,
\log \Gve) =&& -C \{
{(2 \pi)^{-1}} \log |\BGx|^{-1} + \Gz(\BGx) \}
\nonumber \\
&&+ C \{
{(2 \pi)^{-1}} \log \Gve +  H(0,0) \} - H(0,\By),
\label{eq61} \end{eqnarray} with $\Gz$
defined by \eq{eq40aa}.

The condition \eq{eq60} is satisfied provided
\beq
C(\By, \log \Gve) =
\fr{H(0,\By)}{H(0,0) + \fr{1}{2 \pi} \log \Gve - \Gz_\infty }.
\eequ{eq62}
 Combining \eq{eq61}, \eq{eq62}, and \eq{eq54a}, we deduce
\beq
H_\Gve(\Bx, \By) =
-{H(0, \By)}
P_\Gve(\Bx)
+ H(\Bx, \By) +
\tilde{H}_\Gve(\Bx, \By),
\eequ{eq62a} where $\tilde{H}_\Gve$
is the remainder term, such that
\beqa
\GD_x  \tilde{H}_\Gve
(\Bx, \By) = 0 , ~~ \Bx, \By \in \GO_\Gve,
\label{eq64a} \\
\tilde{H}_\Gve
(\Bx, \By) = 0, ~~\Bx \in \prt \GO, ~ \By \in \GO_\Gve,
\label{eq65a} \\
\tilde{H}_\Gve
(\Bx, \By) = H(0,\By)-H(\Bx, \By),
~~\Bx \in \prt F_\Gve, ~ \By \in \GO_\Gve, \label{eq66a}
\end{eqnarray}
where the modulus of the right-hand side in \eq{eq66a} is
estimated by $\mbox{Const} ~\Gve$, uniformly with respect to $\Bx
\in \prt F_\Gve^c$ and $\By \in \GO_\Gve$. The maximum principle
leads to the estimate $ |\tilde H(\Bx, \By)| \leq \mbox{Const} ~
\Gve, $ which is uniform for $\Bx, \By \in \GO_\Gve$.

{\em The approximation of
$h_\Gve$} (see \eq{eq52}--\eq{eq54})
also involves the equilibrium potential $P_\Gve$ from Section \ref{s3a}.
The harmonic function $h_\Gve$ satisfies the homogeneous Dirichlet
condition on $\prt \GO$, and the boundary condition on $\prt
F_\Gve^c$ is rewritten as
$$
h_\Gve(\Bx, \By) = -(2 \pi)^{-1} \log(\Gve^{-1} |\Bx - \By|) - (2
\pi)^{-1} \log \Gve, ~~ \Bx \in \prt F_\Gve^c, \By \in \GO_\Gve.
$$
Hence $h_\Gve(\Bx, \By)$ is sought in the form \beq h_\Gve(\Bx,
\By)= h(\Gve^{-1} \Bx, \Gve^{-1} \By) - (2 \pi)^{-1} \log \Gve +
\tilde{h}_{\Gve}^{(1)}(\Bx, \By), \eequ{seek_h} where the harmonic
function $\tilde{h}_{\Gve}^{(1)}$ vanishes when $\Bx \in \prt
F_\Gve^c, ~ \By \in \GO_\Gve$, and \beq
\tilde{h}_{\Gve}^{(1)}(\Bx, \By) = (2 \pi)^{-1} \log \Gve
-h(\Gve^{-1} \Bx, \Gve^{-1} \By), ~~\Bx \in \prt \GO, \By \in
\GO_\Gve. \eequ{bc_ht}
Representing the right-hand side in \eq{bc_ht} according to Lemma
3, we obtain
$$
\tilde{h}_{\Gve}^{(1)}(\Bx, \By) = (2 \pi)^{-1} \log |\Bx| +
\Gz(\Gve^{-1} \By) + O(\Gve),
$$
uniformly for $\Bx \in \prt \GO, \By \in \GO_\Gve.$ Using the
capacitary potential $P_\Gve$ and the definition \eq{eq3} of
$H(\Bx, \By)$, we
write $\tilde{h}_{\Gve}^{(1)}$ as \beq \tilde{h}_{\Gve}^{(1)}(\Bx,
\By) = - H(\Bx, 0) + \Gz(\Gve^{-1} \By) (1 - P_\Gve(\Bx)) +
\tilde{h}_{\Gve}^{(2)}(\Bx, \By), \eequ{th1} where
$\tilde{h}_{\Gve}^{(2)}$ is  a harmonic function, which is
$O(\Gve)$ for all $\Bx \in \prt \GO, \By \in \GO_\Gve,$ and
satisfies
$$
~~~~~ \tilde{h}_{\Gve}^{(2)}(\Bx, \By) = H(\Bx, 0) = H(0,0) +
O(\Gve),
$$
for all $\Bx \in \prt F_\Gve^c, \By \in \GO_\Gve.$ Hence, \beq
~~~~~~~~~~~~~~~~~~~~~\tilde{h}_{\Gve}^{(2)}(\Bx, \By) = H(0,0)
P_\Gve(\Bx) + O(\Gve), \eequ{th2} uniformly with respect to $\Bx,
\By \in \GO_\Gve$.

Combining \eq{seek_h}, \eq{th1}  and \eq{th2}, we deduce \beqa
h_\Gve(\Bx, \By) =&& h(\Gve^{-1} \Bx, \Gve^{-1} \By) - (2
\pi)^{-1} \log \Gve - H(\Bx, 0) \nonumber \\
&&+ \Gz(\Gve^{-1} \By) (1- P_\Gve(\Bx)) + H(0,0) P_\Gve(\Bx) +
O(\Gve), \label{h_total}
\end{eqnarray}
uniformly with respect to $\Bx, \By \in \GO_\Gve$.

Furthermore, it follows from \eq{eq48a}, \eq{eq62a} and
\eq{h_total} that {Green's function $G_\Gve$} admits the
representation \beqa G_\Gve(\Bx, \By) =&& (2 \pi)^{-1} \log |\Bx -
\By|^{-1} -H(\Bx,
\By) - h(\Gve^{-1} \Bx, \Gve^{-1} \By) \nonumber \\
&& + (2 \pi)^{-1} \log \Gve - \Gz(\BGn) + H(\Bx, 0) \nonumber \\
&& - P_\Gve(\Bx)(H(0,0) - H(0, \By) - \Gz(\Gve^{-1} \By)) +
O(\Gve), \label{G_1}
\end{eqnarray}
which is uniform with respect to $\Bx, \By \in \GO_\Gve$.

By Lemma 4, \eq{G_1} takes the form \beqa \hspace{-.5cm}
G_\Gve(\Bx, \By) &&= (2 \pi)^{-1} \log |\Bx - \By|^{-1} -H(\Bx,
\By) - h(\Gve^{-1} \Bx, \Gve^{-1} \By) \nonumber \\
+ && \fr{(H(0,0) - H(\Bx,0) - \Gz(\Gve^{-1}\Bx))(H(0,0) - H(0,
\By) - \Gz(\Gve^{-1}\By))}{\fr{1}{2 \pi} \log \Gve + H(0,0) -
\Gz_\infty}   \nonumber \\
+ && (2 \pi)^{-1} \log \Gve + H(\Bx,0) + H(0,\By) - H(0,0) +
O(\Gve). \label{G_2}
\end{eqnarray}
Also with the use of Lemma 4, for all $\Bx, \By \in \GO_\Gve$, the
above formula can be written as \beqa G_\Gve(\Bx, \By) = &&(2
\pi)^{-1} \log |\Bx - \By|^{-1} -H(\Bx,
\By) - h(\Gve^{-1} \Bx, \Gve^{-1} \By) \nonumber \\
&&+  ((2 \pi)^{-1} \log \Gve + H(0,0) - \Gz_\infty)
(1-P_\Gve(\Bx))(1-P_\Gve(\By)) \nonumber \\
&&+  (2 \pi)^{-1} \log \Gve + H(\Bx,0) + H(0,\By) - H(0,0)
+ O(\Gve) \nonumber \\
= &&(2 \pi)^{-1} \log |\Bx - \By|^{-1} -H(\Bx,
\By) - h(\Gve^{-1} \Bx, \Gve^{-1} \By) \nonumber \\
&&+  ((2 \pi)^{-1} \log \Gve + H(0,0) - \Gz_\infty)
P_\Gve(\Bx)P_\Gve(\By) \nonumber \\
&& - \Gz(\Gve^{-1} \Bx)- \Gz(\Gve^{-1} \By)+ \Gz_\infty + O(\Gve),
\label{G_3}
\end{eqnarray}
which is equivalent to \eq{eq48}. The proof is complete. $\Box$



\section{Corollaries}
\label{s4cor}

The asymptotic formulae of sections 2 and 3 can be simplified
under
constraints
on positions of the
points $\Bx, \By$ within $\GO_\Gve$.

{\bf Corollary 1.}

{\it (a) Let $\Bx$ and $\By$ be points of $\GO_\Gve \subset {\Bbb
R}^n, n > 2,$ such that \beq ~~~~~~~~~~~~~~~~~\min\{|\Bx|, |\By|\}
> 2 \Gve.\eequ{minxy} Then \beqa G_\Gve(\Bx, \By)  = G(\Bx, \By)
&& - \Gve^{n-2} \mbox{\rm cap}(F) ~
G(\Bx, 0) G(0, \By)
\nonumber \\
&& +O\Big(\fr{\Gve^{n-1}}{(|\Bx| |\By|)^{n-2} \min \{|\Bx|, |\By|
\}}\Big). \label{gcor2a} \end{eqnarray}

(b) If $\max\{|\Bx|, |\By|\} < 1/2$,
then
$$
 G_\Gve(\Bx, \By) = \Gve^{2-n}{g}({\Gve}^{-1}{\Bx}, {\Gve}^{-1}{\By})
$$
\beq - {\cal H}(0,0) (P(\Gve^{-1} \Bx) -1)(P(\Gve^{-1} \By) -1) +
O(\max \{|\Bx|, |\By|  \}).
\eequ{gcor4}
Both {\rm \eq{gcor2a}} and {\rm \eq{gcor4}}
are uniform with respect to $\Gve$ and $(\Bx, \By) \in \GO_\Gve
\times \GO_\Gve$. }


{\bf Proof.}

{\it (a)} The formula \eq{eq16} is equivalent to
\beq
G_\Gve(\Bx, \By) = G(\Bx, \By) - \Gve^{2-n }h(\Gve^{-1} \Bx,
\Gve^{-1} \By)
\eequ{eq16cor} \vspace{-.8cm}
$$
+{\cal H}(0, \By) P(\Gve^{-1} \Bx) + {\cal H}(\Bx, 0) P(\Gve^{-1}
\By)
- {\cal H}(0,0) P(\Gve^{-1} \Bx) P(\Gve^{-1} \By)
$$
$$
-\Gve^{n-2} \mbox{cap} (F) ~{\cal H}(\Bx, 0) {\cal H}(0, \By) +
O\Big(\fr{\Gve^{n-1}}{(\min \{|\Bx|, |\By| \}^{n-2})}\Big).
$$
By Lemmas 1 and 2 \beq P(\Gve^{-1}\Bx) = \fr{\Gve^{n-2} ~\mbox{\rm
cap}(F)}{(n-2) |S^{n-1}| |\Bx|^{n-2}} +
O\Big(\fr{\Gve^{n-1}}{|\Bx|^{n-1}}\Big). \eequ{eq_P} and
\beq \Gve^{2-n} h(\Gve^{-1} \Bx, \Gve^{-1} \By) = \fr{ P(\Gve^{-1}
\By)
}{(n-2) |S^{n-1}| |\Bx|^{n-2}} + O\Big(\fr{
\Gve^{n-1}}{|\Bx|^{n-1} |\By|^{n-2}} \Big) \eequ{eq_h}
\vspace{-.3cm}
$$
=\fr{\Gve^{n-2} \mbox{\rm cap}(F)}{((n-2) |S^{n-1}|)^2 |\Bx|^{n-2}
|\By|^{n-2}} + O\Big(\fr{\Gve^{n-1}}{(|\Bx| |\By|)^{n-2} \min
\{|\Bx|, |\By|\}} \Big).
$$
Direct substitution of \eq{eq_h} and \eq{eq_P} into \eq{eq16cor}
leads to
\beqa G_\Gve(\Bx, \By) =&& G(\Bx, \By) - \fr{\Gve^{n-2}
\mbox{\rm cap}(F)
}{(n-2)^2 |S^{n-1}|^2 |\Bx|^{n-2}
|\By|^{n-2}} \nonumber \\
&&+ \Gve^{n-2} \mbox{\rm cap}(F) \Big(\fr{H(0,\By)}{(n-2)
|S^{n-1}| |\Bx|^{n-2}} + \fr{H(\Bx,0)}{(n-2) |S^{n-1}|
|\By|^{n-2}} \nonumber \\
&& - H(\Bx, 0) H(0, \By) \Big) +
O\Big(\fr{\Gve^{n-1}}{(|\Bx| |\By|)^{n-2} \min \{|\Bx|, |\By|\}}
\Big)
\nonumber \\
=&& G(\Bx, \By) - \Gve^{n-2} \mbox{\rm cap}(F) \Big\{ \Big(
(n-2)^{-1} |S^{n-1}|^{-1} |\Bx|^{2-n} -H(\Bx, 0) \Big) \nonumber
\\
&& \times \Big( (n-2)^{-1} |S^{n-1}|^{-1} |\By|^{2-n} -H(0,
\By)\Big)
\nonumber \\
&& + O\Big(\fr{\Gve^{n-1}}{(|\Bx| |\By|)^{n-2} \min \{|\Bx|,
|\By|\}} \Big),
\nonumber
\end{eqnarray}
which is equivalent to \eq{gcor2a}.

{\it (b)} Since $H(\Bx, \By)$ is smooth in the vicinity of $(\BO,
\BO)$ formula \eq{eq16} can be presented in the form \beqa
G_\Gve(\Bx, \By) =&& \Gve^{2-n} g(\Gve^{-1} \Bx,
\Gve^{-1} \By) - H(0,0) \nonumber \\
&& + (H(0,0) + O(|\By|))P(\Gve^{-1} \Bx) + (H(0,0) + O(|\Bx|))
P(\Gve^{-1} \By) \nonumber \\
&& - H(0,0) P(\Gve^{-1} \Bx) P(\Gve^{-1} \By) + O(\max\{ |\Bx|,
|\By|\}), \nonumber
\end{eqnarray}
which is equivalent to \eq{gcor4}. The proof is complete. $\Box$

We give an analogue of Corollary 1 for the planar case.

{\bf Corollary 2.} {\it
%
(a) Let $\Bx$ and $\By$ be points of $\GO_\Gve \subset {\Bbb R}^2$
subject to \eq{minxy}.
{\it Then
\beq G_\Gve(\Bx, \By) = G(\Bx, \By) +
\fr{G(\Bx,0) G(0, \By)}{\fr{1}{2 \pi} \log \Gve + H(0,0) -
\Gz_\infty} +
 O\Big( \fr{\Gve}{\min \{ |\Bx|,|\By|  \}}\Big),
\eequ{gcor1}

(b) If $\max\{|\Bx|, |\By|\} < 1/2$,
then
\beqa G_\Gve(\Bx, \By) =&& {g}({\Gve}^{-1}{\Bx}, {\Gve}^{-1}{\By})
\nonumber \\
&& + \fr{\Gz(\Gve^{-1}\Bx) \Gz(\Gve^{-1}\By)}{\fr{1}{2 \pi} \log
\Gve + H(0,0) - \Gz_\infty} + O(\max \{|\Bx|, |\By| \}),
\label{gcor3}
\end{eqnarray} Both {\rm \eq{gcor1}} and {\rm \eq{gcor3}} are uniform with
respect to $\Gve$ and $(\Bx, \By) \in \GO_\Gve \times \GO_\Gve $.
} }

{\bf Proof.}
{\it (a)} Formula \eq{eq48} can be written as

$$
G_\Gve(\Bx, \By) = 
({2 \pi})^{-1} \log |\Bx - \By|^{-1} - H(\Bx, \By) - h(\BGx, \BGn)
$$
$$
+ \fr{(G(\Bx, 0) - \Gz(\BGx) + \fr{1}{2 \pi} \log|\BGx| +
\Gz_\infty)(G( 0, \By) - \Gz(\BGn) + \fr{1}{2 \pi} \log|\BGn| +
\Gz_\infty)}{\fr{1}{2 \pi} \log \Gve + H(0,0) - \Gz_\infty}
$$
\vspace{-.4cm} \beq ~~~~~~~~~~~~~~~~~-\Gz(\BGx) - \Gz(\BGn) +
\Gz_\infty +O(\Gve).\eequ{cora}

It follows from Lemma 3 and definition \eq{eq40aa} that \beq
h(\BGx, \BGn) = - (2 \pi)^{-1} \log |\BGx| - \Gz(\BGn) +
O(\Gve/|\Bx|), \eequ{cora1} and \beq \Gz(\BGx) = (2 \pi)^{-1} \log
|\BGx| + \Gz_\infty + O(\Gve/ |\Bx|). \eequ{cora2} Direct
substitution of \eq{cora1} and \eq{cora2} into \eq{cora} yields
\beqa G_\Gve(\Bx, \By) =&& (2 \pi)^{-1} \log |\Bx- \By|^{-1} -
H(\Bx, \By) \nonumber \\
&& + \fr{(-G(\Bx,0) + O(\Gve/|\Bx|))(-G(0,\By)+
O(\Gve/|\By|))}{\fr{1}{2 \pi} \log \Gve + H(0,0) - \zeta_\infty}
+ O(\Gve),
\end{eqnarray}
and hence we arrive at \eq{gcor1}.

{\it (b)} When $\max\{|\Bx|, |\By|\} < 1/2$,
\eq{eq48} is presented in the form: \beqa G_\Gve(\Bx, \By) =&&
g(\Gve^{-1} \Bx, \Gve^{-1} \By) - H(\Bx, \By)
\nonumber \\
\nonumber \\
&& +\fr{(H(0,0) - H(\Bx, 0) - \Gz(\Gve^{-1} \Bx) ) (H(0,0) - H(0,
\By) - \Gz(\Gve^{-1} \By) )}{\fr{1}{2 \pi} \log \Gve + H(0,0) -
\Gz_\infty}  \nonumber \\
\nonumber \\
&& +H(\Bx, 0)+ H(0, \By) - H(0,0) + O(\Gve) \nonumber
\end{eqnarray}
(compare with \eq{G_2}).
Since $H(\Bx, \By)$ is smooth in a vicinity of $(\BO,\BO)$,
we obtain \vspace{-.1cm} \beqa G_\Gve(\Bx, \By) =&& g(\Gve^{-1}
\Bx, \Gve^{-1} \By) + \fr{(- \Gz(\Gve^{-1} \Bx) + O(|\Bx|))(-
\Gz(\Gve^{-1} \By) + O(|\By|))}{\fr{1}{2 \pi} \log \Gve + H(0,0) -
\Gz_\infty} \nonumber \\
&& + O(\max\{|\Bx|, |\By| \})
 \nonumber \\
 \nonumber \\
= &&
g(\Gve^{-1} \Bx, \Gve^{-1} \By) \nonumber \\
&& + \fr{\Gz(\Gve^{-1} \Bx) \Gz(\Gve^{-1} \By) + O(|\By|
\log(|\Bx|/\Gve)) + O(|\Bx| \log(|\By|/\Gve))}{\fr{1}{2 \pi} \log
\Gve + H(0,0) -
\Gz_\infty} \nonumber \\
&& + O(\max\{|\Bx|, |\By| \}), \nonumber
\end{eqnarray}
which implies \eq{gcor3}. $\Box$

\vspace{.1in}



\end{document}

%% file: macroarx.tex
\newcommand\eq[1] {(\ref{#1})}

\newcommand{\bfm}[1]{\mbox{\boldmath ${#1}$}}

\newcommand{\beqa}{\begin{eqnarray}}
\newcommand{\eeqa}[1]{\label{#1}\end{eqnarray}}
\newcommand{\bequ}{\begin{equation}}
\newcommand{\eequ}[1]{\label{#1}\end{equation}}
\newcommand{\Grad}{\nabla}

\newcommand{\ov}[1]{\overline{#1}}

\newcommand{\Gd}{\delta}

\newcommand{\Gve}{\varepsilon}

\newcommand{\Gs}{\sigma}

\newcommand{\Go}{\omega}
\newcommand{\Gx}{\xi}

\newcommand{\Gz}{\zeta}
\newcommand{\GD}{\Delta}

\newcommand{\GO}{\Omega}

\newcommand{\BGn}{\bfm\eta}

\newcommand{\BGx}{\bfm\xi}

\newcommand{\CP}{{\cal P}}

\def\Bl{{\bf l}}

\def\Bn{{\bf n}}

\def\Br{{\bf r}}

\def\Bx{{\bf x}}
\def\By{{\bf y}}
\def\Bz{{\bf z}}
\def\BA{{\bf A}}

\def\BL{{\bf L}}

\def\BO{{\bf O}}

\def\BR{{\bf R}}

\newcommand{\beq}{\begin{equation}}
\newcommand{\eeq}{\end{equation}}
\newcommand{\overliner}{\begin{eqnarray}}
\newcommand{\earr}{\end{eqnarray}}
\newcommand{\beqn}{\begin{equation*}}
\newcommand{\eeqn}{\end{equation*}}
\newcommand{\overlinern}{\begin{eqnarray*}}
\newcommand{\earrn}{\end{eqnarray*}}
\newcommand{\prt}{\partial}

\newcommand{\fr}{\frac}